\documentclass[10pt,a4paper]{amsart}

\usepackage{amssymb,enumerate}
\usepackage{amsmath}

\newcommand{\R}{\mathbb{R}}
\newcommand{\N}{\mathbb{N}}

\newtheorem{theorem}{Theorem}
\newtheorem{lemma}{Lemma}

\newtheorem{proposition}{Proposition}
\newtheorem{definition}{Definition}
\newtheorem{example}{Example}

\begin{document}

\title{Minimal modified energy control\\
for fractional linear control systems\\
with the Caputo derivative}

\thanks{Submitted 17-Dec-2009;
Revised 17-March-2010; recommended 
31-March-2010 for publication 
in \emph{Carpathian J. Math.} 26 (2010), no.~2.}

\author{Dorota Mozyrska}

\address{Dorota Mozyrska\newline
\indent Bia{\l}ystok University of Technology\newline
\indent Faculty of Computer Science \newline
\indent 15-351 Bia\l ystok, Poland}

\email{d.mozyrska@pb.edu.pl}

\author{Delfim F. M. Torres}

\address{Delfim F. M. Torres\newline
\indent University of Aveiro\newline
\indent Department of Mathematics \newline
\indent 3810-193 Aveiro, Portugal}

\email{delfim@ua.pt}




\begin{abstract}
Fractional control systems with the Caputo derivative
are considered. The modified controllability Gramian
and the minimum energy optimal control problem
are investigated. Construction of minimizing steering
controls for the modified energy functional are proposed.
\end{abstract}

\maketitle


\pagestyle{myheadings}
\markboth{Dorota Mozyrska and Delfim F. M. Torres}{Minimal modified energy
control for fractional linear control systems with the Caputo derivative}


\section{Introduction}

The Fractional Calculus is an important Mathematical discipline
\cite{oldman,podlubny,samko}. Several recent books on the
subject have been written, illustrating the usefulness of the theory
in applications \cite{Das,kilbas,S:A:TM,W:B:G,Zaslavsky}.
Different notions of fractional-order derivatives
are available, including the Riemann--Liouville,
the Gr\"{u}nwald--Letnikov, and Caputo,
as well as the generalized functions approach.
One of the youngest fractional derivative,
formulated in 1967, is the Caputo derivative.
The main advantage in using the Caputo derivative is that it avoids problems connected
with initial conditions of fractional differential equations. Indeed, when an
initial value problem is formulated for the Caputo derivative,
initial conditions look like in the classical (non-fractional) case.

The state-space description of fractional order control systems is developed
in \cite{bettayeb,bettayeb:L,kaczorek1,kaczorek2,sierociuk,vinagre}.
Controllability and observability of finite-dimensional fractional differential
systems are investigated mainly in \cite{matignon},
for systems with Riemann--Liouville derivatives.
For positive continuous-time linear systems, the reachability
property has been worked out in \cite{kaczorek1}.
Differently, here we deal with the controllability property
of finite-dimensional linear fractional-order differential systems via the Caputo derivative.
The novelty is the construction of control laws corresponding
to the minimal modified energy, where we use a neutralizer of the singularity
similarly to the one defined in \cite{matignon}.

The paper is organized as follows.
In Section~\ref{sec:II} we recall the main properties of the generalized
Mittag--Leffler function and the fractional analogue for the exponential matrix.
We also review some properties of fractional integrals and derivatives.
In Section~\ref{section:systems} we discuss
the Gramian formula for fractional control systems
and we prove a new result about the minimal modified
energy steering control law (Theorem~\ref{th:1}).
Illustrative examples are discussed. In Section~\ref{sec:rankC:SC}
we recall the classical conditions for controllability
and propose the construction of a different control law
for fractional control systems (Theorem~\ref{th:constr_contr}).


\section{Preliminaries}
\label{sec:II}

We investigate control systems with commensurate order $\alpha \in (0,1]$.
For $\alpha=1$ the classical situation is obtained;
the novelty appears when $\alpha \in (0,1)$.

\begin{definition}[\cite{eld}]\rm
We denote by $E_{\alpha,\beta}$ the two-parameter Mittag--Leffler function
defined by the series expansion
\begin{equation}
\label{eq:ML}
E_{\alpha,\beta}(z)=\sum_{k=0}^{\infty}\frac{z^{k}}{\Gamma(k\alpha +\beta)}\, ,
\quad \alpha>0\, ,  \quad \beta>0 \,.
\end{equation}
When $\beta = 1$, we use the notation $E_{\alpha} = E_{\alpha,1}$.
Let $A\in\R^{n\times n}$.
We extend (\ref{eq:ML}) to the matrix case as follows:
\begin{equation*}
E_{\alpha,\beta}(At^{\alpha})=\sum_{k=0}^{\infty}A^k\frac{t^{k\alpha}}{\Gamma(k\alpha +\beta)}\, .
\end{equation*}
\end{definition}

\begin{definition}[\cite{kilbas}]\rm
Let $A\in \R^{n\times n}$. By
\begin{equation*}
e_{\alpha}^{At}=t^{\alpha-1}\sum_{k=0}^{\infty}A^k\frac{t^{k\alpha}}{\Gamma[(k+1)\alpha]}
= \sum_{k=0}^{\infty}A^k\frac{t^{(k+1)\alpha-1}}{\Gamma[(k+1)\alpha]}
=t^{\alpha-1}E_{\alpha,\alpha}(At^{\alpha})
\end{equation*}
we denote the $\alpha$-exponential matrix function.
\end{definition}

For $\alpha =1$ we have $E_{1}(At)=e_{1}^{At}=\exp(At)$,
where $\exp$ denotes the classical exponential matrix.
We remark that the properties
$\exp(A+B)t = \exp(At)\exp(Bt)$ and
$\exp^{-1}(At) = \exp(-At)$,
satisfied for square matrices $A$ and $B$,
are not valid for functions $e_{\alpha}^{At}$
and $E_{\alpha}(At^{\alpha})$.

\begin{definition}[\cite{kilbas,podlubny,samko}]\rm
Let $\varphi\in L_1\left([t_0,t_1],\R\right)$. The integrals
\begin{gather*}
I^{\alpha}_{t_0+}\varphi(t)
=\frac{1}{\Gamma(\alpha)}\int_{t_0}^t \varphi(\tau)(t-\tau) ^{\alpha-1}d\tau\,, \quad t>t_0,\\
I^{\alpha}_{t_1-}\varphi(t)
=\frac{1}{\Gamma(\alpha)}\int_{t}^{t_1} \varphi(\tau)(\tau-t) ^{\alpha-1}d\tau\,, \quad t< t_1\, ,
\end{gather*}
where $\Gamma$  is  the gamma function and $\alpha>0$,
are called, respectively, the left-sided and the right-sided fractional integrals
of order $\alpha$. Additionally, we define the identity operator
$I\hspace{-1.7mm}I$ by $I\hspace{-1.7mm}I := I^0_{t_0+}=I^0_{t_1-}$.
\end{definition}

We have the following integration by parts formula for fractional integrals.
\begin{proposition}[\cite{kilbas}]
Let $\alpha>0$ and $1/p+1/q\leq 1+\alpha$, $p\geq 1$, $q\geq 1$, with $p\neq 1$
and $q\neq 1$ in the case $1/p+1/q=1+\alpha$.
Then, for $\varphi\in L_p\left([t_0,t_1],\R\right)$
and $\psi\in L_q\left([t_0,t_1],\R\right)$,
the following equality holds:
\begin{equation}
\label{intbypart1}
\int_{t_0}^{t_1}\varphi(\tau)I^{\alpha}_{t_0+}\psi(\tau)d\tau
=\int_{t_0}^{t_1}\psi(\tau)I^{\alpha}_{t_1-}\varphi(\tau)d\tau\,.
\end{equation}
\end{proposition}

\begin{definition}[\cite{kilbas,podlubny}]\rm
Let $\varphi$ be defined on the interval $[t_0,t_1]$.
The left-sided Riemann--Liouville derivative of order $\alpha$
with lower limit $t_0$ is defined by
\begin{equation}
\label{eq:LSRL}
D^{\alpha}_{t_0+}\varphi(t)=\frac{1}{\Gamma(n-\alpha)} \left(\frac{d}{dt}\right)^n
\int_{t_0}^t  \varphi(\tau)(t-\tau)^{n-\alpha-1}d\tau \, ,
\end{equation}
where $n$ is the natural number satisfying $n=[\alpha]+1$
($[\alpha]$ denotes the integer part of  $\alpha$).
Similarly, the right-sided Riemann--Liouville derivative of order $\alpha$
with upper limit $t_1$ is defined by
\begin{equation}
\label{eq:RSRL}
D^{\alpha}_{t_1-}\varphi(t)=\frac{1}{\Gamma(n-\alpha)} \left(-\frac{d}{dt}\right)^n
\int_{t}^{t_1}  \varphi(\tau)(\tau-t)^{n-\alpha-1}d\tau \, .
\end{equation}
\end{definition}

The next proposition is based on \cite[Corollary~2, p.~46]{samko}
and is particularly useful for our purposes.

\begin{proposition}[Integration by parts]
\label{intbypart2}
Let $f\in I^{\alpha}_{t_1-}(L_p)$
and $g\in I^{\alpha}_{t_0+}(L_q)$ with $1/p+1/q\leq 1+\alpha$.
The following formula holds:
\begin{equation}
\label{eq:intbypart}
\int_{t_0}^{t_1}f(t)D^{\alpha}_{t_0+}g(t)dt
=\int_{t_0}^{t_1}g(t)D^{\alpha}_{t_1-}f(t)dt\, , \quad  0<\alpha<1 \, .
\end{equation}
\end{proposition}

\begin{proof}
Denote $D^{\alpha}_{t_1-}f(t)=\varphi(t)$
and $D^{\alpha}_{t_0+}g(t)=\psi(t)$.
The equality (\ref{eq:intbypart}) follows from (\ref{intbypart1})
since $I^{\alpha}_{t_0+}D^{\alpha}_{t_0+}f(t)=f(t)$
is valid for $f\in I^{\alpha}_{t_0+}(L_1)$ (\textrm{cf.} \cite{samko}).
\end{proof}

In the description of fractional control systems we use
the notion of Caputo derivative, which is the preferred
fractional derivative among Engineers.

\begin{definition}[\cite{kilbas}]\rm
The left- and right-sided Caputo fractional derivatives
of order $\alpha\geq 0$ on $[t_0, t_1]$,
denoted by $\mbox{}^CD^{\alpha}_{t_0+}\varphi(t)$
and $\mbox{}^CD^{\alpha}_{t_1-}\varphi(t)$, respectively,
are defined via the Riemann--Liouville
fractional derivatives (\ref{eq:LSRL}) and (\ref{eq:RSRL}) by
\begin{equation}
\label{eq:C}
\mbox{}^CD^{\alpha}_{t_0+}\varphi(t)
= D^{\alpha}_{t_0+}\left(\varphi(t)
-\sum_{k=0}^{n-1}\frac{\varphi^{(k)}(t_0)(t-t_0)^k}{k!}\right)\,
\end{equation}
and
\begin{equation}
\label{eq:C1}
\mbox{}^CD^{\alpha}_{t_1-}\varphi(t)=
D^{\alpha}_{t_1-}\left(\varphi(t)
-\sum_{k=0}^{n-1}\frac{\varphi^{(k)}(t_1)(t_1-t)^k}{k!}\right)\, ,
\end{equation}
where $n=[\alpha]+1$ for $\alpha\not\in \N_0$
and $n=\alpha$ for $\alpha\in \N_0$.
\end{definition}

When $0 <\alpha < 1$ the relations
(\ref{eq:C}) and (\ref{eq:C1}) take the following form:
\begin{equation*}
\mbox{}^CD^{\alpha}_{t_0+}\varphi(t)=
D^{\alpha}_{t_0+}\left(\varphi(t)-\varphi(t_0)\right)\, ,
\quad
\mbox{}^CD^{\alpha}_{t_1-}\varphi(t)=
D^{\alpha}_{t_1-}\left(\varphi(t)-\varphi(t_1)\right)\, .
\end{equation*}

Let $AC[t_0, t_1]$ be the space of functions that are absolutely continuous
on $[t_0, t_1]$ and $AC^n[t_0, t_1]$ denote the space of functions
$\varphi$ that have continuous derivatives up to order $n-1$
on $[t_0, t_1]$ and such that $\varphi^{(n-1)} \in AC[a,b]$.

\begin{proposition}[\cite{kilbas}]
\label{kilbas}
Let $\alpha\geq 0$, $n=[\alpha]+1$ if $\alpha\not\in\N_0$
and $n=\alpha$ if $\alpha\in\N_0$.
If $\varphi\in AC^n[t_0, t_1]$, then the Caputo fractional derivatives
$\mbox{}^CD^{\alpha}_{t_0+}\varphi(t)$ and
$\mbox{}^CD^{\alpha}_{t_1-}\varphi(t)$ exist almost everywhere on $[t_0, t_1]$.
Moreover,
\begin{enumerate}
\item[(a)] If $\alpha \not\in \N_0$, then
\begin{equation*}
\begin{split}
\mbox{}^CD^{\alpha}_{t_0+}\varphi(t)&=\frac{1}{\Gamma(n-\alpha)}\int_{t_0}^t
\varphi^{(n)}(\tau)(t-\tau)^{n-\alpha-1}d\tau
=I_{t_0+}^{n-\alpha}\left(\varphi^{(n)}\right)(t)\, ,\\
\mbox{}^CD^{\alpha}_{t_1-}\varphi(t)
&=\frac{(-1)^n}{\Gamma(n-\alpha)}\int_{t}^{t_1}
\varphi^{(n)}(\tau)(\tau-t)^{n-\alpha-1}d\tau
=(-1)^nI_{t_1-}^{n-\alpha}\left(\varphi^{(n)}\right)(t)\, .
\end{split}
\end{equation*}
\item[(b)] If $\alpha=n\in\N_0$, then
$\mbox{}^CD^{n}_{t_0+}\varphi(t)=\varphi^{(n)}(t)$
and $\mbox{}^CD^{n}_{t_1-}\varphi(t)=(-1)^n\varphi^{(n)}(t)$.
\end{enumerate}
\end{proposition}

For a function $x:[0,T]\rightarrow\R^n$ we use
similar notation as in the classical case:
\[
\mbox{}^CD^{\alpha}_{0+}x(t)=\mbox{}^CD^{\alpha}_{0+}
\left(\begin{array}{c}
x_1(t)\\
\vdots \\
x_n(t)
\end{array}\right)
= \left(\begin{array}{c}
\mbox{}^CD^{\alpha}_{0+}x_1(t)\\
\vdots \\
\mbox{}^CD^{\alpha}_{0+}x_n(t)
\end{array}\right) \,.
\]
Such situation, when for each component we use the same
fractional order $\alpha$ of differentiation
(in the Riemann--Liouville or Caputo sense),
is known in the literature as the fractional derivative
with commensurate order \cite{bettayeb,kaczorek1,sierociuk}.

\begin{proposition}
\label{prop:sam}
For $\alpha>0$ the following holds:
\begin{enumerate}
\item[(i)] $\mbox{}^CD^{\alpha}_{t_0+}E_{\alpha}(A(t-t_0)^{\alpha})
=AE_{\alpha}(A(t-t_0)^{\alpha})$;
\item[(ii)] $D^{\alpha}_{t_0+} e_{\alpha}^{A(t-t_0)}
=Ae_{\alpha}^{A(t-t_0)}$;
\item[(iii)] $D^{\alpha}_{T-} S(T-\tau)=AS(T-\tau)$,
where $S(t)=e_{\alpha}^{At}$.
\end{enumerate}
\end{proposition}

\begin{proof}
(i) Directly from the definition of the classical
Mittag--Leffler function $E_{\alpha}$,
and from the formula of the Caputo derivative of a power function, one has:
\[\mbox{}^CD^{\alpha}_{t_0+}(t-t_0)^{\beta}
=\frac{\Gamma(\beta+1)}{\Gamma(\beta-\alpha+1)}(t-t_0)^{\beta-\alpha},
\quad \beta\neq 0\,.
\]
Since the Caputo derivative of a constant function is zero,
\begin{equation*}
\begin{split}
\mbox{}^CD^{\alpha}_{t_0+}E_{\alpha}(A(t-t_0)^{\alpha})
&= \mbox{}^CD^{\alpha}_{t_0+}\sum_{k=0}^{\infty}A^k\frac{(t-t_0)^{k\alpha}}{\Gamma(k\alpha +1)}
=\sum_{k=1}^{\infty}A^k\frac{(t-t_0)^{(k-1)\alpha}}{\Gamma((k-1)\alpha+1)}\\
&=AE_{\alpha}(A(t-t_0)^{\alpha})\,.
\end{split}
\end{equation*}

(ii) As the formula for the $\alpha$-exponential function does not consist
of constant terms, and since the formula for the Riemann--Liouville
derivative of a power function is the same as for the Caputo derivative, we have:
\[
D^{\alpha}_{t_0+} e_{\alpha}^{A(t-t_0)}
=D^{\alpha}_{t_0+}\sum_{k=0}^{\infty}A^k \frac{(t-t_0)^{(k+1)\alpha-1}}{\Gamma[(k+1)\alpha]}
=\sum_{k=1}^{\infty}A^k \frac{(t-t_0)^{k\alpha-1}}{\Gamma(k\alpha)}
=Ae_{\alpha}^{A(t-t_0)}\, ,
\]
where we used the fact that
$\lim\limits_{\alpha\rightarrow 0}\frac{1}{\Gamma(\alpha)}=0$.

(iii) Using the formulas
\[
D^{\alpha}_{T-}(T-\tau)^{\beta-1}=\frac{\Gamma(\beta)}{\Gamma(\beta-\alpha)}(T-\tau)^{\beta-\alpha-1}\, ,
\quad
\lim_{\beta\rightarrow \alpha}D^{\beta}_{T-}\frac{(T-\tau)^{\alpha-1}}{\Gamma(\alpha)}=0\,
\]
(\textrm{cf.}, \textrm{e.g.}, \cite{samko}), we get:
\begin{equation*}
\begin{split}
D^{\alpha}_{T-} S(T-\tau)
&= D^{\alpha}_{T-}\left(I\frac{1}{\Gamma(\alpha)}(T-\tau)^{\alpha-1}
+A\frac{(T-\tau)^{2\alpha-1}}{\Gamma{2\alpha}}+\cdots\right)\\
&= A\frac{(T-\tau)^{\alpha-1}}{\Gamma(\alpha)}+A^2\frac{(T-\tau)^{2\alpha-1}}{\Gamma(2\alpha)}+\cdots\\
&= (T-\tau)^{\alpha-1}A\sum_{k=0}^{\infty} A^k\frac{(T-\tau)^{k\alpha}}{\Gamma[(k+1)\alpha]}\\
&=AS(T-\tau)\, .
\end{split}
\end{equation*}
\end{proof}

\begin{proposition}
For $\alpha>0$ the following relation holds:
\begin{equation*}
E_{\alpha}(A(t-t_0)^{\alpha})
= I+\int_{t_0}^t Ae_{\alpha}^{A(t-\tau)}d\tau\, .
\end{equation*}
\end{proposition}

\begin{proof}
Follows by direct calculation of the integral:
\begin{equation*}
\begin{split}
\int_{t_0}^t Ae_{\alpha}^{A(t-\tau)}d\tau
&= \int_{t_0}^t \sum_{k=0}^{\infty}A^{k+1}\frac{(t
-\tau)^{(k+1)\alpha-1}}{\Gamma[(k+1)\alpha]}d\tau\\
&= \sum_{k=1}^{\infty}A^{k}\frac{(t-t_0)^{k\alpha}}{\Gamma(k\alpha+1)}
=E_{\alpha}(A(t-t_0)^{\alpha})-I\, .
\end{split}
\end{equation*}
\end{proof}

Since $e_{\alpha}^{At}=t^{\alpha-1}E_{\alpha,\alpha}(At^{\alpha})$
and each Mittag--Leffler function $E_{\alpha,\alpha}(az^{\alpha})$,
$\alpha>0$, is an entire function on the complex plane,
we can state the following:
\begin{proposition}
\label{lemma:3}
Let $\alpha>0$.
There is a uniquely determined function
$g(t)=t^{1-\alpha}G(t)$ such that
$e_{\alpha}^{At}g(t)=E_{\alpha,\alpha}(At^{\alpha})G(t)=I$,
for $t\neq 0$, and $\lim\limits_{t\rightarrow0}e_{\alpha}^{At}g(t)=I$.
\end{proposition}

\begin{lemma}
\label{lemma:intbypart}
Let $\alpha>0$ and $\psi(t) \in I^{\alpha}_{0+}(L_q)$,
where  $1/p+1/q\leq 1+\alpha$ for $p$ such that
all components of $S(T-t)$ belong to $I^{\alpha}_{T-}(L_p)$. Then,
\begin{equation*}
\int_0^TS(T-\tau)D^{\alpha}_{0+}\psi(\tau)d\tau
=\int_0^TAS(T-\tau)\psi(\tau)d\tau\,.
\end{equation*}
\end{lemma}

\begin{proof}
Taking into account (\ref{eq:intbypart})
and item (iii) of Proposition~\ref{prop:sam}, we have
\[
\int_0^TS(T-\tau)D^{\alpha}_{0+}\psi(\tau)d\tau
=\int_0^T\psi(\tau)D^{\alpha}_{T-}S(T-\tau)d\tau
= \int_0^TAS(T-\tau)\psi(\tau)d\tau\,.
\]
\end{proof}


\section{Control systems and the Gramian}
\label{section:systems}

We consider the following linear time-invariant control system
of order $\alpha\in(0,1]$, denoted by $\Sigma$:
\begin{equation*}
\mbox{}^CD^{\alpha}_{0+}x(t) = Ax(t)+Bu(t) \, ,
\quad y = Cx(t)\, ,
\end{equation*}
where $x(t)\in \mathbb{R}^n$, $u(t)\in \mathbb{R}^m$,
matrix $A\in M_{n\times n}(\mathbb{R})$,
$B\in M_{n\times m}(\mathbb{R})$,
$C\in M_{p\times n}(\mathbb{R})$,
and $\mbox{}^CD^{\alpha}_{0+}$ indicates the fractional
Caputo derivative of commensurate order $\alpha$.

The forward trajectory of the system $\Sigma$,
starting at $t_0=0$ and evaluated at $t\geq 0$,
is the solution of the initial value problem
$\mbox{}^CD^{\alpha}_{0+} x(t) = Ax(t)+Bu(t)$,
$x(0)=a\in\R^n$ \cite{kilbas}:
\begin{equation}
\label{eq:solution}
\gamma(t,a,u)=\left(I+\int_{0}^tS(\tau) Ad\tau\right)a
+\int_{0}^t S(t-\tau)Bu(\tau)d\tau \, ,
\end{equation}
where $S(t)= e_{\alpha}^{At}$.
Moreover, we can represent (\ref{eq:solution}) in the following way:
\begin{equation*}
\gamma(t,a,u)=S_0(t)a+\int_{0}^t S(t-\tau)Bu(\tau)d\tau \, ,
\end{equation*}
where $S_0(t)=E_{\alpha}(At^{\alpha})=I+\int_{0}^tS(\tau) Ad\tau$.
The formula for the forward trajectory can be obtained
using the Laplace transform \cite{kaczorek1}.
Taking into account the output of $\Sigma$, the forward output trajectory
is then defined by values evaluated at $t\geq  0$:
\[
\eta(t,a,u) = C\gamma(t,a,u)
=C\left(I+\int_{0}^tS(\tau) Ad\tau\right)a
+C\int_{0}^t S(t-\tau)Bu(\tau)d\tau\, .
\]

For system $\Sigma$ we define the notion of controllability
in the standard manner:

\begin{definition}\rm
Let $T>0$. The system $\Sigma$ is controllable on $[0,T]$
if for any $a\in\R^n$ and $b\in\R^n$
there is a control $u(\cdot)$ defined on $[0,T]$
which steers the initial state $\gamma(0,a,u)=a$
to the final state $\gamma(T,a,u)=b$.
\end{definition}

Following \cite{matignon,zabczyk} we denote by
\begin{equation*}
Q_T =\int_{0}^{T}  S(T-t)BB^* S^*(T-t)(T-t)^{2(1-\alpha)}dt
\end{equation*}
the \emph{controllability Gramian}
of fractional order $\alpha$,
on the time interval $[0,T]$,
corresponding to the system $\Sigma$.
As in the classical case \cite{zabczyk}, $Q_{T}$ is symmetric
and nonnegative definite. The term $(T-t)^{2(1-\alpha)}$ under
the integral is called in \cite{matignon} a neutralizer
of the singularity at $t=T$. It is needed in order
to ensure the convergence of the integral.

Let $T>0$. By $L^2_{\alpha}\left([0,T],\R^m\right)$ we denote
the set of functions $\varphi:[0,T]\rightarrow \R^m$ such that
$\widetilde{\varphi}$ defined by
$\widetilde{\varphi}(t)=(T-t)^{\alpha-1}\varphi(t)$
is square integrable on $[0,T]$.

\begin{theorem}
\label{th:1}
Let $T>0$ and $Q_T$ be nonsingular. Then,
\begin{itemize}
\item[(a)] for any states $a$, $b\in \mathbb{R}^n$
the control law
\begin{equation}
\label{eq1}
\overline{u}(t)=-(T-t)^{2(1-\alpha)}B^*S^*(T-t)Q_T^{-1}f_T(a,b) \, ,
\quad t\in[0,T),
\end{equation}
where
\[f_T(a,b)=\left(I+\int_{0}^{T}S(t)Adt\right)a-b=-b+S_0(T)a\]
and $\overline{u}(T)=0$,
drives point $a$ to point $b$ in time $T$;
\item[(b)] among all possible controls from $L^2_{\alpha}\left([0,T],\R^m\right)$
driving $a$ to $b$ in time $T$,
the control $\overline{u}$ defined by~(\ref{eq1})
minimizes the integral
\begin{equation}
\label{eq:}
\int_0^T|(T-t)^{\alpha-1}u(t)| ^2dt \, .
\end{equation}
Moreover,
\[
\int_0^T|(T-t)^{\alpha-1}\overline{u}(t)|^2 dt
= <Q_T^{-1}f_T(a,b), f_T(a,b)>\, ,
\]
\end{itemize}
where $<\cdot,\cdot>$ denotes the inner product.
\end{theorem}

\begin{proof}
(a) From the expression of $\overline{u}$ we have that
\begin{equation*}
\begin{split}
\gamma(T,a,\overline{u}) &= f_T(a,b)+b \\
& \qquad -\left(\int_0^T (T-t)^{2(1-\alpha)}S(T-t)BB^*S^*(T-t)dt\right)Q_T^{-1}f_T(a,b)\\
&= f_T(a,b)+b-Q_TQ_T^{-1}f_T(a,b)=b\,.
\end{split}
\end{equation*}
Moreover,
$\lim\limits_{t\rightarrow T-}\overline{u}(t)=\overline{u}(T)$ since
$\lim\limits_{t\rightarrow T-}(T-t)^{2(1-\alpha)}B^*S^*(T-t)Q_T^{-1}f_T(a,b)$.

(b) We notice that for $h(t)=|(T-t)|^{2(1-\alpha)}$
\begin{equation*}
\begin{split}
\int_0^T & |(T-t)^{\alpha-1}\overline{u}(t)|^2 dt
= \int_0^T\left| (T-t)^{1-\alpha}B^*S^*(T-t)Q_T^{-1}f_T(a,b)\right|^2dt\\
&= \int_0^T h(t)\left<B^*S^*(T-t)Q_T^{-1}f_T(a,b),  B^*S^*(T-t)Q_T^{-1}f_T(a,b)\right>dt\\
&=\left<\int_0^T h(t)S(T-t)BB^*S^*(T-t)dt, Q_T^{-1}f_T(a,b)\right>\\
&=\left<Q_TQ_T^{-1}f_T(a,b),Q_T^{-1}f_T(a,b)\right>
=\left<f_T(a,b),Q_T^{-1}f_T(a,b)\right>\,.
\end{split}
\end{equation*}
Let us take another control $u$ for which
$(T-t)^{\alpha-1}u(t)$ is square integrable
on $[0,T]$ and $\gamma(T,a,u)=b$. Then,
\begin{equation*}
\begin{split}
\int_0^T & (T-t)^{2(\alpha-1)}\left<u(t),\overline{u}(t)\right>dt\\
&=-\int_0^T h(t) \left<u(t),(T-t)^{2(1-\alpha)}B^*S^*(T-t)Q_T^{-1}f_T(a,b)\right>dt\\
&= -\int_0^T\left<u(t),B^*S^*(T-t)Q_T^{-1}f_T(a,b)\right>dt
= \left<f_T(a,b),Q_T^{-1}f_T(a,b)\right>\,.
\end{split}
\end{equation*}
Hence,
\begin{equation*}
\int_0^T(T-t)^{2(1-\alpha)}\left<u(t),\overline{u}(t)\right>dt
=\int_0^T(T-t)^{2(1-\alpha)}\left<\overline{u}(t),\overline{u}(t)\right>dt
\end{equation*}
and from that we obtain
\begin{multline*}
\int_0^T(T-t)^{2(1-\alpha)}|u(t)|^2dt \\
=\int_0^T(T-t)^{2(1-\alpha)}|\overline{u}(t)|^2 dt
+\int_0^T(T-t)^{2(1-\alpha)}|u(t)-\overline{u}(t)|^2 dt \, ,
\end{multline*}
which gives the minimality property for the integral.
\end{proof}

As in the classical case, the operator
$$u\mapsto \gamma(T,0,u)=\int_0^TS(t)Bu(T-t)dt$$
is a linear operator from the space
$L_{\alpha}^2\left([0,T],\R^m\right)$ into $\R^n$.
Hence the classical result can
be replaced by the following proposition.
\begin{proposition}[\textrm{cf.} \cite{zabczyk}]
\label{prop:2}
If any state $b\in\R^n$ is attainable from $a=0$,
then the matrix $Q_T$ is nonsingular for any arbitrary $T>0$.
\end{proposition}

\begin{example}\rm
Let $\Sigma$ be the following system evaluated on $\R^2$:
\[
\Sigma: \quad
\begin{cases}
\mbox{}^CD^{0.5}_{0+}x_1(t)=x_2(t)\, , \\
\mbox{}^CD^{0.5}_{0+}x_2(t)=u(t)\,.
\end{cases}
\]
Let us take $a=(1,0)^*$ and $b=(0,0)^*$.
Since $A=\left(\begin{array}{cc} 0 & 1\\ 0 & 0\end{array}\right)$
and $B=\left(\begin{array}{c} 0 \\  1 \end{array}\right)$,
we obtain that
$S(t)=\left(\begin{array}{cc}\frac{1}{\sqrt{\pi t}} & 1
\\ 0 & \frac{1}{\sqrt{\pi t}}\end{array}\right)$ while
$S_0(t)=\left(
\begin{array}{cc} 1 & \frac{2\sqrt{t}}{\sqrt{\pi}} \\
0 & 1
\end{array}\right)$.
Hence, the formula for the solution with the initial condition $\gamma(0,a,u)=a$ is
\begin{equation*}
\gamma(t,a,u)=\left(
\begin{array}{cc} 1 & \frac{2\sqrt{t}}{\sqrt{\pi}} \\
0 & 1 \end{array}\right)a
+ \int_0^t \left(\begin{array}{cc} \frac{1}{\sqrt{\pi(t-\tau)}} & 1 \\
0 & \frac{1}{\sqrt{\pi(t-\tau)}}
\end{array}\right)Bu(\tau)d\tau\, .
\end{equation*}
Let us take $u(t)\equiv 1$. Then, for the given $a$,
$\gamma(t,a,u)= \left( 1+t \quad 2\frac{\sqrt{t}}{\sqrt{\pi}}\right)^*$.
From the last expression we see that using constant
$u(\cdot)\equiv 1$ for $t>0$
we are not able to steer the given initial
point $a$ to the origin.

Let now $f_T(a,b)=S_0(T)a-b=a$.
The Gramian has the form
\[Q_T=\left(\begin{array}{cc}
\frac{T^2}{2} & \frac{2T^{3/2}}{3\sqrt{\pi}} \\
\frac{2T^{3/2}}{3\sqrt{\pi}} &\frac{T}{\pi}\end{array}\right)\]
and the control \[\overline{u}(t)=-\frac{18(T-t)}{T^2}
+\frac{12\sqrt{T-t}}{T^{3/2}}\] drives $a$ to $b$ with the modified energy
\[m=\int_0^T|(T-t)^{-0.5}\overline{u}(t)|^2dt=\frac{18}{T^2}\,.\]
\end{example}

\begin{example}\rm
Let $\alpha\in(0,1)$.
Consider the following fractional system $\Sigma$ on $\R^2$:
\[\Sigma: \quad
\begin{cases}
\mbox{}^CD^{\alpha}_{0+}x_1(t)=x_2(t), \\
\mbox{}^CD^{\alpha}_{0+}x_2(t)=-x_1(t)+u(t)\,.
\end{cases}
\]
The matrix $A$ is now skew-symmetric, and thus
\[A^0=I, \quad
A=\left(
\begin{array}{rr} 0& 1\\ -1 & 0
\end{array}\right),
\quad A^2=-I\,.
\]
Hence, $A^k = I$ if  $k=0,4,8,\ldots$;
$A^k = A$ if $k=1,5,9,\ldots$;
$A^k = -I$ if $k=2,6,10,\ldots$;
and $A^k =-A$ if $k=3,7,11,\ldots$. Moreover,
\begin{equation*}
\begin{split}
S(t)&=t^{\alpha-1}\left(I\frac{1}{\Gamma(\alpha)}+A\frac{t^{\alpha}}{\Gamma(2\alpha)}-
I\frac{t^{2\alpha}}{\Gamma(3\alpha)}-A\frac{t^{3\alpha}}{\Gamma(4\alpha)}+\cdots\right)\\
&=I\left(\frac{t^{\alpha-1}}{\Gamma(\alpha)}-\frac{t^{3\alpha-1}}{\Gamma(3\alpha)}+\right)
+A\left(\frac{t^{2\alpha-1}}{\Gamma(2\alpha)}-\frac{t^{4\alpha-1}}{\Gamma(4\alpha)}+\cdots\right)\,.
\end{split}
\end{equation*}
Using the notation
\begin{equation*}
\sin_{\alpha}t=\sum_{k=0}^{\infty}(-1)^k\frac{t^{2(k+1)\alpha-1}}{\Gamma[2(k+1)\alpha]}\, ,
\quad
\cos_{\alpha}t=\sum_{k=0}^{\infty}(-1)^k\frac{t^{(2k+1)\alpha-1}}{\Gamma[(2k+1)\alpha]}\, ,
\end{equation*}
we can write $S(t)=\left(\begin{array}{rr}\cos_{\alpha}t
& \sin_{\alpha}t \\ -\sin_{\alpha} t & \cos_{\alpha}t
\end{array}\right)$.
As $B=\left(\begin{array}{c}0 \\ 1\end{array}\right)$,
we have
\begin{equation*}
Q_T = \int_0^T (T-t)^{2(1-\alpha)} \cdot M_\alpha(t) dt
\end{equation*}
with
\begin{equation*}
M_\alpha(t) =
\left(\begin{array}{cc}
\sin^2_{\alpha}(T-t) & \sin_{\alpha}(T-t)\cos_{\alpha}(T-t)\\
\sin_{\alpha}(T-t)\cos_{\alpha}(T-t) & \cos^2_{\alpha}(T-t)
\end{array}\right) \, .
\end{equation*}
To get an exact formula for $Q_T$ and $Q_T^{-1}$
is difficult. We can, however,
easily obtain approximations
with the desired precision
for both matrices,
and an approximate formula for the optimal control.
Let us consider a concrete situation. Let $\alpha=\frac{1}{2}$.
Then $\sin_{\frac{1}{2}}t=e^{-t}$ and
\[\cos_{\frac{1}{2}}t=\frac{1}{\sqrt{\pi t}}\left(1
-\sum_{k=1}^{+\infty}\frac{2^kt^{2k-1}}{\prod_{i=1}^k(2i-1)}\right)
=\frac{1}{\sqrt{\pi t}}(1-2t+\cdots) \,. \]
We choose to approximate $\cos_{\frac{1}{2}}t$ by functions
$$
c_L(t)=\frac{1}{\sqrt{\pi t}}\left(1
-\sum_{k=1}^{L}\frac{2^kt^{2k-1}}{\prod_{i=1}^k(2i-1)}\right) \, ,
\quad L=1, 2, \ldots
$$
In this way we can approximate $Q_T$ by
\begin{equation*}
\int_0^T (T-t)\left(
\begin{array}{cc}
e^{-2(T-t)} & c_L(T-t)e^{-(T-t)} \\
c_L(T-t)e^{-(T-t)} & c_L^2(t)
\end{array}\right) dt\,.
\end{equation*}
Let $T=10$, the starting point to be
$a=\left( 0 \quad 1\right)^*$,
and the final point to be the origin, \textrm{i.e.},
$b= \left(0 \quad 0\right)^*$.
Applying formula (\ref{eq1}) for the steering control $\overline{u}$
with the approximate expressions for $S^*(T-t)$ and $Q_T^{-1}$,
we check that we are very close
to the final goal $b$.
Taking more terms in the functions $c_L$,
the value of the modified energy tends to the minimum.
For example, for $L=1$ we have
$c_1(t)=\frac{1}{\sqrt{\pi t}}(1-2t)$ and the approximate value of energy is
$m_{1}=1.02$; for $L=11$ we have $m_{11}=0.0921$;
and for $L=12$ we get $m_{12}=0.0911$. These numbers were calculated using
MATLAB with a final precision of four digits.
\end{example}


\section{Rank conditions and steering controls}
\label{sec:rankC:SC}

In the classical theory it is quite easy to explain why the condition
$\mbox{rank}\, B=n$ is sufficient for controllability.
Roughly speaking, one needs to construct a
special control $u(\cdot)$. In the case of fractional order systems
we use the special function given by Proposition~\ref{lemma:3}.
\begin{proposition}
\label{prop:control}
Let $\mbox{rank}\, B=n$ and $B^+$ be such that $BB^+=I$.
Let $g(\cdot)$ be the matrix function defined by Proposition ~\ref{lemma:3}.
Then the control
$$
\widehat{u}(t)=\frac{1}{T}B^+g(T-t)\left(b-S_0(T)a\right)\, ,
\quad t\in[0,T] \, ,
$$
transfers $a$ to $b$ in time $T> 0$.
\end{proposition}

\begin{proof}
The proof follows by a direct calculation:
\begin{equation*}
\gamma(T,a,\widehat{u}) = S_0(T)a
+\frac{1}{T}\int_0^T S(T-t)BB^+
g(T-t)\left(b-S_0(T)a\right)dt=b\,.
\end{equation*}
\end{proof}

\begin{example}\rm
Let $\Sigma$ be the following fractional system evaluated on $\R$:
\[\Sigma: \quad \mbox{}^CD^{\alpha}_{0+}x(t)=u(t)\, , \]
where $\alpha\in(0,1]$.
Let us take $a,b\in\R$, $T>0$. Then,
\begin{equation*}
S(t)=\frac{t^{\alpha-1}}{\Gamma(\alpha)}\,,
\quad S_0(t)=1, \quad B=B^+=1 \, ,
\end{equation*}
and
\[\gamma(T,a,u)= a+\frac{1}{\Gamma(\alpha)}\int_0^T (T-t)^{\alpha-1}u(t)dt\,.\]
Consider, accordingly to Proposition~\ref{lemma:3},
$g(t)=t^{1-\alpha}\Gamma(\alpha)$. From Proposition~\ref{prop:control},
\[\widehat{u}(t)=\frac{\Gamma(\alpha)}{T}(T-t)^{1-\alpha}(b-a)\]
transfers $a$ to $b$: $\gamma(T,a,\widehat{u})=b$
with the modified energy~(\ref{eq:}) of value
\[m=\frac{\Gamma^2(\alpha)(b-a)^2}{T}\,.\]
It is also the minimum energy as $\widehat{u}(t)=\overline{u}(t)$.
\end{example}

An algebraic condition equivalent to controllability
for fractional linear control systems has been derived
in \cite{fliess} and cited again in \cite{bettayeb,matignon}.
According to Theorem~\ref{th:1} and  Proposition~\ref{prop:2},
we can state a similar formulation as in \cite{zabczyk}:
\begin{theorem}
The following conditions are equivalent:
\begin{enumerate}
\item[(a)] An arbitrary state $b\in\R^n$ is attainable from $0$.
\item[(b)] System $\Sigma$ is controllable.
\item[(c)] System $\Sigma$ is controllable at a given time $T>0$.
\item[(d)] Matrix $Q_T$ is nonsingular for some $T>0$.
\item[(e)] Matrix $Q_T$ is nonsingular for an arbitrary $T>0$.
\item[(f)] $\mbox{rank}\, \left[ A|B\right]=\mbox{rank}\,
\left[B, AB,\ldots, A^{n-1}B\right]=n$.
\end{enumerate}
\end{theorem}

If the rank condition is satisfied, then the control $\overline{u}(\cdot)$
given by (\ref{eq1}) steers $a$ to $b$ at time $T$.
Our goal now is to find another formula for the steering control
by using the matrix $\left[A|B\right]$ instead of the controllability matrix $Q_T$.
It is a classical result that if $\mbox{rank}\,\left[A|B\right]=n$,
then there exists a matrix $K\in M(mn,n)$ such that $\left[A|B\right]K=I\in M(n,n)$ or,
equivalently, there are matrices $K_1$, $K_2$, $\ldots$, $K_n \in M(m,n)$ such that
$BK_1+ABK_2+\cdots +A^{n-1}BK_n=I$.

For the next construction it is more convenient
to use the notion of Riemann--Liouville derivative.
We begin by introducing a notation for compositions
of the Riemann--Liouville derivatives
with the same order $\alpha$, $\alpha\in(0,1)$.
Let $R^{\alpha,0}_{0+}\psi(t)=\psi(t)$.
Then for $j\in\N$, recursively, we put
$R^{\alpha,j+1}_{0+}\psi(t)
:=D^{\alpha}_{0+}\left(R^{\alpha,j}_{0+}\psi(t)\right)$.

\begin{theorem}
\label{th:constr_contr}
Let $\mbox{rank}\,\left[A|B\right]=n$ and $\alpha\in(0,1)$.
Let $p$ be such that $S(T-t)\in I^{\alpha}_{T-}(L_p)$ and
$\varphi$ be a real function given on $[0,T]$ such that
\begin{enumerate}
\item[(i)] $\int_0^T \varphi(t)dt=1$;
\item[(ii)] $R^{\alpha,j}_{0+}\psi(t)\in I^{\alpha}_{0+}(L_q)$
for $j=0,\ldots, n-1$, where
$\psi(t)=g(t)\left(b-S_0(T)a\right)\varphi(t)$
and $S(T-t)g(t)=I$, $t\in [0,T]$,
for $1/p+1/q\leq 1+\alpha$.
\end{enumerate}
Then the control
$\hat{u}(t)=K_1\psi(t)+K_2D^{\alpha}_{0+}\psi(t)
+\cdots+K_nR^{\alpha, n-1}_{0+}\psi(t)$, $t\in[0,T]$,
transfers a to b at time $T\geq 0$.
\end{theorem}

\begin{proof}
Using $j-1$ times the formula (\ref{eq:intbypart})
of integration by parts and Lemma~\ref{lemma:intbypart},
\begin{equation*}
\int_0^T S(T-t)BK_jR^{\alpha, j-1}_{0+}\psi(t)dt=\int_0^T S(T-t)A^{j-1}BK_j\psi(t)dt
\end{equation*}
for $j=1,\ldots, n$. Then,
\[\int_0^TS(T-t)B\hat{u}(t)dt=\int_0^TS(T-t)B\psi(t)dt\]
and
$\gamma(T,a,\hat{u}) = S_0(T)a+\int_0^TS(T-t)g(t)\left(b-S_0(T)a\right)\varphi(t)dt = b$.
\end{proof}


\section{Conclusion}

The concept of fractional (\textrm{i.e.}, non-integer)
derivative and integral is increasingly being recognized
as a good tool to model the behavior of complex systems
in various fields of science and engineering
\cite{A:M:T:10,AML:Ric:Del,B:08,bhn,El-Nabulsi:Torres,S:A:TM,vinagre}.
This is particularly true in the areas of control theory
and control engineering, where fractional order controllers
and plants provide confirmed evidence of better
performances when compared with the best integer order controllers
\cite{A:B:07,B:D:A:09,bettayeb,Das,gastao:delfim,bettayeb:L,kaczorek1,kaczorek2,matignon}.
In this work we consider linear time-invariant control systems
of fractional order $\alpha$, $\alpha \in  (0,1)$, in the sense
of Caputo. The controllability question
for such systems is addressed.
Main results give sufficient conditions
assuring the existence of controls that
drive the system between any two desired states.
Explicit formulas for such controls are given.
A control law is obtained
that minimizes the modified energy
$\int_0^T|(T-t)^{\alpha-1}u(t)|^2 dt$.
Our results are reduced to the
classical ones when the fractional order
of differentiation $\alpha$ tends to one.


\subsection*{Acknowledgments}

Dorota Mozyrska was supported by Bia{\l}ystok University of Technology grant S/WI/1/08;
Delfim F. M. Torres by the R\&D unit ``Centre for Research on Optimization and Control'' (CEOC)
of the University of Aveiro, cofinanced by the European Community Fund FEDER/POCI 2010.

We are very grateful to two anonymous referees for valuable suggestions.



\end{document}